\documentclass[12pt]{article}
\usepackage{preview}
\usepackage{amsmath}
\usepackage{mathdots}
\usepackage{amsfonts}
\usepackage{amssymb}
\usepackage{latexsym}
\usepackage{graphicx}
\usepackage{pgfplots}
\usepackage{todonotes}
\usepackage{tikz}
\usetikzlibrary{tikzmark,calc}
\newcommand*{\Scale}[2][4]{\scalebox{#1}{\ensuremath{#2}}}%
\usepackage{color}
\usepackage{fullpage}
       \newcommand{\rowsize}{.7cm}
       \newcommand{\colsize}{2.5cm}

\def\eig{\mbox{eig}_{\lambda\neq 0}}

\newtheorem{prop}{Proposition}

\newcommand{\rank}{\mbox{rank}}

\newcommand{\ignore}[1]{}

\title{
The low-rank eigenvalue problem
}
\author{Yuji Nakatsukasa\thanks{{\tt nakatsukasa@maths.ox.ac.uk}, Mathematical Institute, University of Oxford, Oxford, OX2 6GG, UK, and National Institute of Informatics.}
}

\begin{document}
\maketitle

\begin{abstract}
The nonzero eigenvalues of $AB$ are equal to those of $BA$: 
an identity that holds as long as the products are square, 
 even when $A,B$ are rectangular. 
This fact naturally suggests an efficient algorithm for computing eigenvalues and eigenvectors 
of a low-rank matrix $X= AB$ with $A,B^T\in\mathbb{C}^{N\times r}, N\gg r$: form the small $r\times r$ matrix $BA$ and find its eigenvalues and eigenvectors. 
For nonzero eigenvalues, the eigenvectors are related by $ ABv = \lambda v \Leftrightarrow  BAw = \lambda w $ with $w=Bv$, and the same holds for Jordan vectors. 
For zero eigenvalues, the Jordan blocks can change sizes between $AB$ and $BA$, and we characterize this behavior. 

\end{abstract}

\section{Introduction}
Low-rank matrices are omnipresent in scientific computing, 
often due to the need of compressing data and the fact that they allow for efficient algorithms. The literature on low-rank matrices is too vast to list. 

 This note concerns the efficient computation of the eigenvalues and eigenvectors of low-rank matrices. The literature on this subject appears to be surprisingly sparse and incomplete; for example the discussion in~\cite{venit2000eigenvalues} is based on characteristic polynomials and not computationally easy to use. 

\section{Algorithm}
Suppose $X$ is a large but low-rank matrix such that 
$X= AB$, where $A\in\mathbb{C}^{N\times r}$, $B\in\mathbb{C}^{r\times N}$ with $N\gg r$.
The key identity we rely on is $\eig(AB)=\eig(BA)$, 
where $\eig(X)$ denotes the set of nonzero eigenvalues of $X$, counting multiplicities. 
This is a  classical result that can be proved as follows
(a proof attributed to Kahan in~\cite[p.~27]{stewart-sun:1990}): Let 
$Y =
\begin{bmatrix}
I& 0\\-B &I 
\end{bmatrix}
$. Then we have 
\begin{equation}
  \label{eq:Kahan}
Y^{-1}\begin{bmatrix} AB& A\\0& 0_{r\times r}  \end{bmatrix}Y=
\begin{bmatrix}
I&0 \\B &I 
\end{bmatrix}
\begin{bmatrix}
 AB& A\\0& 0_{r\times r} 
\end{bmatrix}
\begin{bmatrix}
I&0 \\-B& I 
\end{bmatrix}
=\begin{bmatrix}
0_{N\times N} &A \\0 & BA
\end{bmatrix},   
\end{equation} 
so the matrices 
$\big[\begin{smallmatrix}AB& A\\0& 0_{r\times r}   \end{smallmatrix}\big]$ and 
$
\big[
\begin{smallmatrix}
0_{N\times N} &A \\0 & BA  
\end{smallmatrix}
\big]
$ are similar, thus proving the identity. 

Clearly, converting from $AB$ to $BA$ may be advantageous since the eigenvalues of $BA$ can be computed with $O(r^3)$ flops as opposed to $O(N^3)$ (we also need $2Nr^2$ flops to compute $BA$). We illustrate this as follows: 
\[
\mbox{eig}_{\lambda\neq 0} \left(\ 
   \begin{tikzpicture}[node distance=1mm,baseline=-1mm]
       \node[draw, rectangle, solid, fill=gray!10, minimum width=\rowsize, minimum height=\colsize] (A) {$\Scale[1]{A}$};
       \node[draw, rectangle, solid, fill=gray!10, minimum width=\colsize, minimum height=\rowsize,right =of A] (B) {$\Scale[1]{B}$};
   \end{tikzpicture}
\right)=
\mbox{eig}_{\lambda\neq 0}\left(\ 
   \begin{tikzpicture}[node distance=1mm,baseline=-1mm]
       \node[draw, rectangle, solid, fill=gray!10, minimum width=\colsize, minimum height=\rowsize] (B) {$\Scale[1]{B}$};
       \node[draw, rectangle, solid, fill=gray!10, minimum width=\rowsize, minimum height=\colsize,right =of B] (A) {$\Scale[1]{A}$};
   \end{tikzpicture}
\ \right)=\mbox{eig}_{\lambda\neq 0}\left(\ 
   \begin{tikzpicture}[node distance=1mm,baseline=-1mm]
       \node[draw, rectangle, solid, fill=gray!10, minimum width=\rowsize, minimum height=\rowsize] (B) {$\Scale[.8]{BA}$};
   \end{tikzpicture}
\ \right).
\]

In addition to the eigenvalues, the eigenvectors and indeed Jordan chains of $AB$ can also be obtained easily from those of $BA$ (which is almost immediate from~\eqref{eq:Kahan}). 




\ignore{
In addition to $\mbox{eig}(AB)=\mbox{eig}(BA)$ for nonzero eigenvalues, this proof establishes further that 
(i) the size of the Jordan blocks are all the same for $\lambda\neq 0$, and 
(ii) the eigenvectors are related as follows. Let $ABx = \lambda x$ with $A$ tall. Then 
$BA(Bx) = \lambda (Bx)$, so $Bx$ is an eigenvector of $BA$. 
The opposite also holds: 
if $BAy = \lambda y$, then 
$AB(Ay) = \lambda (Ay)$, so $Ay$ is an eigenvector of $AB$. 
}

\begin{prop}\label{ABeignonzero}
The nonzero eigenvalues of $AB\in\mathbb{C}^{N\times N}$ and $BA\in\mathbb{C}^{r\times r}$ are identical. 
If $(\lambda,v)$ is an eigenpair of $BA$ with $\lambda\neq 0$, then 
$(\lambda,Av)$ is an eigenpair of $AB$. 
Moreover, if $v_1,\ldots,v_k$ forms a Jordan chain for $AB$ with eigenvalue $\lambda$, then 
$Av_1,\ldots,Av_k$ is a Jordan chain for $AB$ for the same eigenvalue. 
\end{prop}
{\sc proof.} 
The equality of eigenvalues was established above in~\eqref{eq:Kahan}. 

Suppose $BAv = \lambda v$ with $v\neq 0$, $\lambda\neq 0$. 
Then $ABAv = \lambda Av$, so defining $w=Av$, we have $ABw = \lambda w$. 
Furthermore, $w\neq 0$, because $w=0$ implies $BAv = 0$, contradicting $\lambda\neq 0$. 

The same argument holds for Jordan chains: 
Suppose that $V_1=[v_1,\ldots,v_k]$ forms a Jordan chain for $BA$, so that 
$BAV_1=V_1J$, where $J$ is a Jordan block with eigenvalue $\lambda\neq 0$. 
Then $(AB)AV_1=AV_1J$, so 
$AV_1=[Av_1,\ldots,Av_k]$ (which has full column rank $k$) forms a Jordan chain for $AB$. 
\hfill $\square$

We note that since the statement and argument are symmetric about $A$ and $B$ and we did not assume $N\geq r$, 
the same proof shows that if $ABw = \lambda w$ with $w\neq 0$, $\lambda\neq 0$, then  $BAv = \lambda v$ with $v=Bw$, and likewise for Jordan chains.

Proposition~\ref{ABeignonzero} immediately suggests an efficient algorithm for computing eigenvalues, eigenvectors and Jordan chains
of a low-rank matrix $X\in\mathbb{C}^{N\times N}$:
\begin{enumerate}
\item Find $A,B$ such that $X= AB$ (or $X\approx AB$), with $A,B^T\in\mathbb{C}^{N\times r}$ with $r\leq N$. 
\item Compute the matrix product $BA$, and its nonzero eigenvalues $\{\lambda_i\}_{i=1}^{r_0}$, and if desired, 
eigendecomposition 
$BAV=V\Lambda$ where $\Lambda=\mbox{diag}(\lambda_1,\ldots,\lambda_{r_0})$, $V=[v_1,\ldots,v_{r_0}]\in\mathbb{C}^{r\times r_0}$, or Jordan decomposition $BAV=VJ$. 
\item The nonzero eigenvalues and eigenvectors of $AB$ are $(\lambda_i,Av_i)$, with eigenvalue decomposition $ABW=W\Lambda$ 
(where $W=[w_1,\ldots,w_{r_0}]\in\mathbb{C}^{N\times r_0}$ with $w_i=Av_i$)
or
Jordan decomposition  $ABW=WJ$. 
\end{enumerate}

How to perform the first step is outside the scope of this note; a large number of algorithms have been proposed. A popular approach is the randomized SVD~\cite{halko2011finding}. 
Of course, when $X$ is given in low-rank form as a product $X=AB$, one can skip the first step. 

Once $A,B$ are available, the cost of computing eigenvalues is  $2Nr^2+cr^3$ flops, where $2Nr^2$ is for forming $BA$ and the scalar $c$ is about $\frac{4}{3}$ if $BA$ is symmetric or $9$ otherwise, and larger (respectively $9$ and $\geq 25$) when eigenvectors are required~\cite{golubbook4th}. 
It is almost trivial that the above approach is highly efficient. To begin with, a naive approach would require $cN^3$ flops (though flops are merely one measure of the algorithm complexity, not necessarily reflecting the running time). 
Compared with other solvers such as Arnoldi, the cost is lower and more predictable: Arnoldi costs $O(Nr^2+r^3)$ if one uses the low-rank structure for matrix-vector multiplications, but with a larger constant with no known iteration count for convergence; mathematically $r$ steps of Arnoldi would suffice. The above approach is also clearly simpler than Arnoldi, both to understand and implement. 

\subsection{Preserving symmetry}
When $X=AB$ is symmetric (or Hermitian, for which replace $^T$ with $^*$ below), it is desirable to preserve the symmetry in the algorithm, in particular the realness of the eigenvalues and orthogonality of the eigenvectors. 

A symmetric low-rank matrix has a decomposition $X=\tilde A\tilde S\tilde A^T$ where $\tilde S\in\mathbb{R}^{r\times r}$ is symmetric. Such decomposition may be computed for example as in~\cite[\S~5.1]{martinsson2016randomized}
By a congruence transformation $S=W\tilde SW^T$ for a nonsingular $W\in\mathbb{R}^{r\times r}$, we can reduce $S$ to a diagonal matrix of $\pm 1$'s, and set $X=ASA^T$ with $A:=\tilde AW^{-1}$. 

We would then compute the eigenvalues of $A^TAS$---but this is nonsymmetric. One can work around this by 
noting that $(A^TAS,I)S=(A^TA,S)$, and solving the equivalent $r\times r$ \emph{generalized} eigenvalue problem $A^TAv=\lambda Sv$, which 
 is symmetric positive definite assuming $\rank(A)=r$; if not, we have $\mbox{rank}(X)=\tilde r<r$, and a lower-rank representation $X=\tilde A\tilde B$ with $\tilde A,\tilde B^T\in\mathbb{C}^{N\times \tilde r}$ is possible. We then compute an eigendecomposition $A^TAV=SV\Lambda$---exploiting symmetry, for example via the Cholesky factorization of $A^TA$~\cite[Ch.~8]{golubbook4th}---such that $V^T(A^TA,S)V=(I_r,\Lambda^{-1})$, noting that $\Lambda$ is invertible because $A^TA$ is positive definite. Now left-multiplying $AS$ gives $ASA^TAV=ASSV\Lambda=AV\Lambda$, so taking $W:=AV$ we have $ASA^TW=W\Lambda$. Note that $W$ has orthonormal columns, since $W^TW=V^TA^TAV=I_r$. We thus obtain the nonzero eigendecomposition for $X=ASA^T=W\Lambda W^T$. 

\section{Zero eigenvalues}
Above we have fully covered the nonzero eigenvalues of $AB$ and $BA$; in practice, this is most likely all that matters in applications of low-rank matrices, as the zero eigenvalue has high multiplicity, at least $N-r$. 
Nonetheless, a mathematically interesting question is to relate the zero eigenvalues and eigenvectors of $AB$ and $BA$ (if any). This leads to a classical result  by Flanders~\cite{flanders1951elementary}, who discovered the interesting phenomenon that the Jordan chains of $AB$ and $BA$ for zero eigenvalues can have lengths differing by 1, unlike the situation for nonzero eigenvalues; this has been extended to products of more matrices in~\cite{de2014flanders}. 

We examine the behavior of zero eigenvalues of $AB$ and $BA$ in our low-rank context. We assume that $A,B^T$ are both of full column-rank; otherwise we can further reduce $r$ as discussed above. 

\begin{prop}\label{prop0}
Let $AB\in\mathbb{C}^{N\times N}$ and $BA\in\mathbb{C}^{r\times r}$ ($N\geq r$), with $\mbox{rank}(A)=\mbox{rank}(B)=r$. Suppose that 0 is an eigenvalue of $BA$ with geometric multiplicity $\ell$, with Jordan blocks of size $k_1,\ldots, k_\ell$. Then 
$AB$ has eigenvalue 0 with geometric multiplicity $N-r+\ell$, with 
Jordan block sizes 1 ($N-r-\ell$ copies), and $k_1+1,\ldots, k_\ell+1$. 
\end{prop}
{\sc proof}. 
Suppose that $BAV=VJ$ with $J$ a $k\times k$ Jordan block with eigenvalue 0, and $\mbox{rank}(V)=k$. 
By assumption $A$ is of full column rank, so we have  $\mbox{rank}(AV)=k$, 
 and left-multiplying by $A$ gives $AB(AV)=(AV)J$; this is a valid Jordan chain for $AB$. The question is, is $k$ the full length of this Jordan chain? We claim that the answer is no. To see this, first note that 
$\mbox{rank}(BA)=r-\ell$ by the definition of $\ell$. 
Since $\mbox{rank}(AB)=r$ by assumption, 
and the nonzero eigenvalues do not change their Jordan block sizes between $AB$ and $BA$, it follows that $\mbox{rank}(AB)$ must be equal to $r-\ell$ plus the increase in the Jordan block sizes for $\lambda=0$ from $BA$ to $AB$. 
Now by Flanders' theorem, 
zero eigenvalues can change their Jordan block size by at most one, so it follows that each of the $\ell$ Jordan blocks of $BA$ (including $1\times 1$ blocks) must increase its length by one, noting that a $1\times 1$ Jordan block $(=0)$ of $AB$ does not increase the rank. 


These, along with the nonzero eigenvalues, account for dimension $r+\ell$. Since the matrix $AB$ has a null space of dimension $N-r$, of which $\ell$ are of the form $Av$ where $BAv=0$, there is a subspace of dimension $N-(r+\ell)$ remaining in the null space; these account for the remaining zero eigenvalues of $AB$. 
\hfill$\square$

To illustrate the rank increase, consider the case where $A=a,B=b^T$ are vectors with $b^Ta=0$. Then $ab^Ta=0$ but $ab^Tb=\kappa a$ for some nonzero scalar $\kappa$, so $[a,b]$ is a Jordan chain for the matrix $ab^T$. 
 It is perhaps interesting that in our low-rank setting, a zero eigenvalue of the small matrix $BA$ \emph{always} implies a growth in the Jordan block size of $AB$. As the rank-one example $AB=ab^T$ suggests, unlike eigenvalues the singular values of $AB$ and $BA$ are not related, zero or nonzero. 

We note that the lengths of the Jordan blocks described in Propositions~\ref{ABeignonzero} (for $\lambda\neq 0$) and~\ref{prop0} (for $\lambda= 0$)
 sum up to the full dimension $N$; we have thus 
identified the complete Jordan structure of $AB$. 

Let us close with a remark on classical algorithms for computing (not necessarily low-rank) eigenvalues, namely the QR and LR algorithms~\cite{wilkinson:1965}. In their simplest form, they are based on decomposing the matrix $X=AB$ and swapping the order $\tilde X=BA$, and repeating (usually $O(N)$ times, employing acceleration techniques including shifts). The algorithm we presented for low-rank matrices does the same!---but with nonsquare factors, and taking just one step.





\subsection*{Acknowledgment}
I thank Abi Gopal and Nick Trefethen for their comments on a draft. 

\bibliographystyle{abbrv}
\bibliography{C:/Users/Nakatsukasa/Dropbox/paper/bib2}

\end{document}